\documentclass{article}
\usepackage{amsmath, amsthm, amsfonts, amssymb}
\usepackage{hyperref}
\usepackage[latin1]{inputenc}

\makeatletter \@addtoreset{equation}{section}

\makeatletter

\def\beg   {\begin{theorem}}   \def\ee   {\end{theorem}}
\def\be   {\begin{equation}}   \def\ee   {\end{equation}}
\def\ba   {\begin{array}}      \def\ea   {\end{array}}
\def\bea  {\begin{eqnarray}}   \def\eea  {\end{eqnarray}}
\def\bean {\begin{eqnarray*}}  \def\eean {\end{eqnarray*}}

\newtheorem{lemma}{Lemma}[section]
\newtheorem{theorem} [lemma]{Theorem}

\begin{document}

\vspace{4cm}
\begin{center} \LARGE{\textbf{ Uniqueness of Some Differential Polynomials of Meromorphic Functions }}
 \end{center}
 \vspace{1cm}
 \begin{center} \bf{Kuldeep Singh Charak$^{1 }$, \quad Banarsi Lal$^{2}$ }
\end{center}

\begin{center}
 Department of Mathematics, University of Jammu,
Jammu-180 006, INDIA.\\
{$^{1}$ E-mail: kscharak7@rediffmail.com }\\

{$^{2}$ E-mail: banarsiverma644@gmail.com }
\end{center}

\bigskip
\begin{abstract}
In this paper, we prove some uniqueness results  which improve and generalize several earlier works. Also, we prove a value distribution result concerning $f^{(k)}$ which provides a partial answer to a question of Fang and Wang [A note on the conjectures of Hayman, Mues and Gol'dberg, \emph{Comp. Methods, Funct. Theory} (2013)  \textbf{13}, 533--543].
\end{abstract}

\vspace{1cm}\noindent \textbf{Keywords: } Meromorphic functions, small functions, sharing of values, Nevanlinna theory.

\vspace{0.5cm} \noindent\textbf{AMS subject classification: 30D35, 30D30}

\vspace{4cm}

\normalsize
\newpage

\section{Introduction} Throughout, by a meromorphic function we always mean a non-constant meromorphic function in the complex plane $\mathbb C$. We use the standard notations of Nevanlinna Theory such as $m(r,f),~N(r,f),~T(r,f),~S(r,f)$ etc. (one may refer to \cite{HAY}). Let $f$ and $g$ be two meromorphic functions and $a \in \mathbb C$. By $E(a,f)$, we denote the set of zeros of $f - a$ counting multiplicities (CM) and by $\overline E(a, f)$, the set of zeros of $f - a$ ignoring multiplicities (IM). Two meromorphic functions $f$ and $g$ are said to share the value $a$ CM if $E(a,f) = E(a,g)$ and to share the value $a$ IM if $\overline E(a,f) = \overline E(a,g)$. Further, by $E_{k)}(a,f)$, we denote the set of zeros of $f-a$ with multiplicities atmost $k$ in which each zero is counted according to its multiplicity. Also, by $\overline E_{k)}(a,f)$, we denote the set of zeros of $f - a$ with multiplicity atmost $k$, counted once.

We denote by $\mathcal{A}$, the class of meromorphic functions $f$ satisfying $$\overline{N}(r,f) + \overline{N}(r, \frac{1}{f}) = S(r,f).$$  Clearly, each member of class $\mathcal{A}$ is a transcendental meromorphic function. Further, by $\mathcal M(D)$ we denote the space of all meromorphic functions on a domain $D$. A mapping $M: \mathcal M(\mathbb C) \rightarrow \mathcal M(\mathbb C)$ given by $$M[f] = a.\prod^{k}_{j = 0}(f^{(j)})^{n_j};\;\forall\; f \in \mathcal M(\mathbb C)$$ with $n_0,\;n_1,\dotsc, n_k$ as non-negative integers and $a \in \mathcal M(\mathbb C)\;:\;a \not\equiv 0$; is called a {\it differential monomial} of {\it degree} $d = \sum^{k}_{j = 0} n_j$ and {\it weight} $w(M) = \sum^{k}_{j = 1} (1 + n_j)$. We call $a$ the co-efficient of $M$. If $a = 1$, then $M$ is said to be {\it normalised}.  A sum $P: = \sum^{p}_{j = 1}M_j$ of differential monomials $M_1,\;M_2,\dotsc, M_p$ which are linearly independent over $\mathcal M(\mathbb C)$ is called a {\it differential polynomial} of degree $$\text{deg}(P)\;:= \text{max}\{\text{deg}(M_1),\;\text{deg}(M_2),\dotsc,\text{deg}(M_p)\}$$ and the weight $$w(P)\;:\;=\text{max}\{w(M_1),\;w(M_2),\dotsc w(M_p)\}.$$ If deg$(M_1) = \dotsb =$ deg $(M_p) = d$, we call $P$, {\it homogeneous} (of degree $d$).  Also for any $a \in \mathbb C$, we define
$$N_1\left(r,\frac{1}{f - a}\right) = N\left(r,\frac{1}{f - a}\right) - \overline N\left(r,\frac{1}{f - a}\right).$$
and
$$N_2\left(r, \frac{1}{f-a}\right) = \overline N\left(r,\frac{1}{f-a}\right) + \overline N_{(2}\left(r, \frac{1}{f-a}\right),$$
where $N_{(k}\left(r, 1/(f - a)\right)$ is the counting function of those zeros of $f - a$ whose multiplicity  is atleast $k$, and $\overline N_{(k}\left(r,1/(f - a)\right)$ is the one corresponding to ignoring multiplicity.
Finally, by $S(f)$, we denote the set of small functions of $f$; that is,
$$S(f): = \{a \;|\; a \;\text {is meromorphic and} \;T(r,a) = S(r,f)\; \text {as}\; r \rightarrow \; \infty \}.$$

The uniqueness theory of meromorphic functions has perfected the value distribution theory of Nevanlinna and has a vast range of applications in Complex Analysis. Particularly, uniqueness theory of meromorphic functions has been proved to be a handy tool in dealing with the problems on normal families of meromorphic functions. For recent progress concerning normality, one may refer to  \cite{QD}, \cite{XX}, and \cite{YXL}. For recent developments in the uniqueness theory of meromorphic functions (sharing, weighted sharing and q-difference sharing of polynomials), one may refer to \cite{LY}, \cite{Y},\cite{YL} and \cite{ZY}.\\

In the present paper, we prove some uniqueness results which improve and generalize the works of Yang and Yi \cite{CM} , Wang and Gao \cite{GAO}, and Huang and Huang \cite{HHB}. Also, a partial answer to a question of Fang and Wang \cite{FM} concerning value distribution of $f^{(k)} - a$, where $k \in \mathbb N$ and $a \;(\not\equiv 0, \infty)$ is a small function of $f$, is obtained.

\section{Main Results}

 Yang and Yi \cite[Theorem 3.29, p.197]{CM} proved the following result for class  $\mathcal{A} $:

\begin{theorem}
Let $f,g \in \mathcal{A} $, and $a$ be a non-zero complex number. Furthermore, let $k$ be a positive integer.\\
$(i)$ If $\overline E_{1)}(a, f) = \overline E_{1)}(a, g)$ , then $f \equiv g$ or $f.g \equiv a^2$.\\
$(ii)$ If $\overline E_{1)}(a, f^{(k)}) = \overline E_{1)}(a, g^{(k)})$, then $f \equiv g$ or $f^{(k)}.g^{(k)} \equiv  a^2$.
\end{theorem}

A function $f$ is said to \emph{share a value $a$ partially with $g$ IM} if $\overline E(a,f) \subseteq \overline E(a,g)$. We use the notation $N_{1)}\left(r, 1/(g-a)| f\neq a\right)$, to denote the simple zeros of $f-a$, that are not the zeros of $g-a$. Using this notation and the notion of partial sharing,   we improve Theorem $2.1$ as:

\begin{theorem}
Let $f , g \in \mathcal{A} $, $a$ be a non-zero complex number and $k$ be a positive integer.\\
$(i)$ If $\overline E_{1)}(a, f) \subseteq \overline E_{1)}(a, g)$  and $N_{1)}\left(r, 1/(g-a)| f\neq a\right) = S(r,g)$, then $f \equiv g$ or $f.g \equiv a^2$.\\
$(ii)$ If $\overline E_{1)}(a, f^{(k)}) \subseteq \overline E_{1)}(a, g^{(k)})$ and $N_{1)}\left(r, 1/(g^{(k)}-a)| f^{(k)}\neq a\right) = S(r,g)$, then $f \equiv g$ or $f^{(k)}.g^{(k)} \equiv  a^2$.
\end{theorem}

{\bf Example}. Consider $f(z)=e^z$ and $g(z)=e^{2z}$. Then $f, g \in \mathcal{A}$, $\overline E_{1)}(1, f) \subseteq \overline E_{1)}(1, g)$ and $N_{1)}(r, 1/(g-1)| f\neq 1) \neq S(r,g)$, and the conclusion of Theorem $2.2$ does not hold. Thus, the condition ``$N_{1)}\left(r, 1/(g-a)| f\neq a\right) = S(r,g)$" in Theorem $2.2$, is essential.

Theorem $2.2$ also holds if we take $a$ to be a small function different from $0$ and $\infty$, as in that case we can take functions $F = f/a$ and $G = g/a$ instead of $f$ and $g$ so that $F, G \in \mathcal{A}$.\\

  In $2011$, Huang and Huang \cite[Theorem 3, p. 231]{HHB} improved a result of Yang and Hua \cite[Theorem 1, p. 396]{CH} as
\begin{theorem}
Let $f$ and $g$ be two meromorphic functions and $n \geq 19$ be an integer. If $E_1(1, f^nf') = E_1(1,g^ng')$, then either $f=dg$ for some $(n+1)$th root of unity $d$ or  $f(z)=c_1e^{cz}$ and $g(z)=c_2e^{-cz}$, where $c,c_1,c_2$ are constants satisfying $(c_1c_2)^{(n+1)}c^2 = -1$.
\end{theorem}

 In this paper, we improve Theorem $2.3$ for functions of  class $ \mathcal{A} $ as
\begin{theorem}
Let $f, g \in \mathcal{A}$, $n \geq 2$ be an integer and $a (\neq0) \in \mathbb C$. If $\overline E_{1)}(a, f^nf') = \overline E_{1)}(a, g^ng')$, then either $f=dg$ for some $(n+1)$th root of unity $d$ or  $f(z)=c_1e^{cz}$ and $g(z)=c_2e^{-cz}$, where $c,c_1,c_2$ are constants satisfying $(c_1c_2)^{(n+1)}c^2 = -a^2$.
\end{theorem}

Concerning sharing of small functions, Wang and Gao \cite[Theorem 1.3, p.2]{GAO} proved
\begin{theorem}
Let $f$ and $g$ be two transcendental meromorphic functions, $a (\not\equiv 0) \in S(f) \cap S(g)$, and let $n \geq 11$ be positive integer. If $f^nf'$ and $g^ng'$ share $a$ CM, then either $f^nf'g^ng' \equiv a^2$, or $f = tg$ for a constant such that $t^{n+1} =1$.
\end{theorem}

Here in this paper, we partially extend this result to a more general class of differential polynomials as
\begin{theorem}
Let $f$ and $g$ be two transcendental meromorphic functions, $a (\not\equiv 0) \in S(f) \cap S(g)$, and let $n, m, k$ be positive integers satisfying $n > km + 3m + 2k + 8$, and $m > k-1$. If $f^n(f^m)^{(k)}$ and $g^n(g^m)^{(k)}$ share $a$ CM, then either
$$f^n(f^m)^{(k)}g^n(g^m)^{(k)} \equiv a^2\; \text{or}\; f^n(f^m)^{(k)} \equiv g^n(g^m)^{(k)}.$$
\end{theorem}
For $m > k-1$, we have $n > k^2 + 4k + 5$ so that by substituting $k = 1$, we get $n >10$. Thus Theorem $2.6$ reduces to Theorem $2.5$.\\

For the differential polynomials, Barker and Singh \cite[Theorem 3, p.190] {baker} proved
\begin{theorem}
The differential equation
$$af^nf' + P_{n -1}(f) = 0,$$
where $a (\not\equiv 0) \in S(f) $ has no transcendental meromorphic solution $f$ satisfying $N(r,f) = S(r,f)$, where $P_{n -1}(f)$ is a homogeneous differential polynomial of degree $n-1$.
\end{theorem}
In a similar way, we can prove the following more general result
\begin{theorem}
The differential equation
$$af^n(f^m)^{(k)} + P_{n -1}(f) = 0,$$
where $a (\not\equiv 0) \in S(f)$ and $m,\;n$ are positive integers, has no transcendental meromorphic solution $f$ satisfying $N(r,f) = S(r,f)$, where $P_{n -1}(f)$ is a homogeneous differential polynomial of degree $n-1$.
\end{theorem}

 Concerning the value distribution of $k$th derivative of a meromorphic function, Fang and Wang \cite[Proposition 3, p.542]{FM} proved the following result:

\begin{theorem}
Let $f$ be a transcendental meromorphic function having atmost finitely many simple zeros. Then $f^{(k)}$ takes on every non-zero polynomial infinitely often for $k = 1,\;2,\;3,\;\dotsc .$.
\end{theorem}
In the same paper Fang and Wang \cite[Question 2, p.543]{FM} posed the following question:

{\bf Question}: {\it Let $f$ be a transcendental meromorphic function having atmost finitely many simple zeros. Must $f^{(k)}$ take on every non-zero rational function infinitely often for $k = 1,\;2,\;3,\;\dotsc .$. ?}\\

Here, we give a partial answer to this question involving small function as
\begin{theorem}
Let $f$ be a transcendental meromorphic function having atmost finitely many simple zeros and $N \left( r, 1/f''\right) = S(r, f)$. Let $a (\not\equiv 0, \infty) \in S(f)$, then $f^{(k)} - a$ has infinitely many zeros for $k = 1,\;2,\;3,\;\dotsc .$.
\end{theorem}

\section{Some Lemmas}
 We recall the following results which we shall use in the proof of main results of this paper:
\begin{lemma} \label{111} \cite[Theorem 3, p.396] {CH}
Let $f$ and $g$ be two non-constant entire functions, $n\geq 1$ and $a (\neq 0) \in \mathbb C$. If $f^nf'g^ng' = a^2$, then $f(z)=c_1e^{cz}$ and $g(z)=c_2e^{-cz}$, where $c,c_1,c_2$ are constants satisfying $(c_1c_2)^{(n+1)}c^2 = -a^2$.
\end{lemma}
 \begin{lemma} \label{0} \cite[Lemma 1.10, p.82]{CM}
 Let $f_1$ and $f_2$ be non-constant meromorphic functions and $c_1, c_2 \;\text{and}\;c_3$ be non-zero constants. If $c_1f_1 + c_2f_2 \equiv c_3$, then
 $$ T(r, f_1) < \overline N\left(r, \frac{1}{f_1}\right) + \overline N\left(r, \frac{1}{f_2}\right) + \overline N(r, f_1) + S(r,f_1).$$
 \end{lemma}

\begin{lemma}  \label{22} \cite[Lemma 3.8, p.193] {CM}
If $f \in \mathcal{A}$ and $k$ is a positive integer, then $f^{(k)} \in \mathcal{A}$.
\end{lemma}

\begin{lemma} \label{33} \cite [Lemma 3.9, p.194]{CM}
If $f, g \in \mathcal{A}$ and $f^{(k)} = g^{(k)}$, where $k$ is a positive integer, then $f \equiv g$.
\end{lemma}

\begin{lemma} \label{44} \cite [Lemma 3.10, p.194]{CM}
If $f\in \mathcal{A}$ and $a$ is a finite non-zero number, then
$$N_{1)}\left(r, \frac{1}{f - a}\right) = T(r, f) + S(r,f),$$
where $N_{1)}\left(r, 1/(f - a)\right)$ denotes the simple zeros of $f - a$.
\end{lemma}

\begin{lemma} \label{55} \cite[Theorem 1.24, p.39] {CM}
Suppose $f$ is a nonconstant meromorphic function and $k$ is a positive integer. Then
$$N\left(r, \frac{1}{f^{(k)}}\right) \leq N\left(r, \frac{1}{f}\right) + k \overline N(r, f) + S(r, f).$$
\end{lemma}

\begin{lemma} \label{66} \cite [Lemma 2.3, p.3]{GAO}
Let $f$ and $g$ be two meromorphic functions . If $f$ and $g$ share 1 CM, then one of the following must occur:\\
(i) $ T(r, f) + T(r, g) \leq 2\{N_2\left(r, 1/f\right) + N_2\left(r, 1/g\right) + N_2(r,  f) + N_2(r, g)\} + S(r,f) + S(r, g)$,\\
(ii) either $f\equiv g$ or $fg\equiv1$.
\end{lemma}

\begin{lemma} \cite[Lemma 1, p.537]{FM}\label{68}
 Let $f$ be a transcendental meromorphic function, $k \geq 2$ be an integer, and $\epsilon > 0$. Then
 $$(k - 1)\overline N(r, f) + N_1\left(r, \frac{1}{f}\right) \leq N\left(r, \frac{1}{f^{(k)}}\right) + \epsilon T(r, f).$$
 \end{lemma}

\section{Proof of Main Results}

\begin{proof}[Proof of Theorem 2.2]
Since $\overline E_{1)}(a, f) \subseteq \overline E_{1)}(a, g)$,
$$N_{1)}\left(r, \frac{1}{f-a}\right) \leq N_{1)}\left(r, \frac{1}{g-a}\right).$$
Since (by Lemma \ref {44})
$$N_{1)}\left(r, \frac{1}{f-a}\right) = T(r, f) + S(r,f)$$
and
$$N_{1)}\left(r, \frac{1}{g-a}\right) = T(r, g) + S(r,g),$$
therefore,
$$N_{(2}\left(r, \frac{1}{f - a}\right) = S(r, f),$$
$$ N_{(2}\left(r, \frac{1}{g - a}\right) = S(r, g)$$
 and
\begin{equation} \label{a}
T(r,g) \geq T(r, f) + S(r,f).
\end{equation}
Define a function $h : \mathbb C \rightarrow \overline {\mathbb C}$ by
\begin{equation} \label{b}
h(z) = \frac{f(z)-a}{g(z)-a}.
\end{equation}
Since $\overline E_{1)}(a, f) \subseteq \overline E_{1)}(a, g)$, we have
\begin{equation} \label{c}
\overline N(r, h) \leq \overline N(r, f) + \overline N_{(2}\left(r, \frac{1}{g-a}\right) + N_{1)}\left(r, \frac{1}{g-a}| f \neq a\right) = S(r,g)
\end{equation}
\begin{equation} \label{d}
\overline N(r, \frac{1}{h}) \leq \overline N(r, g) + \overline N_{(2}\left(r, \frac{1}{f-a}\right) = S(r,g)
\end{equation}
and
$$T(r,h) \leq T(r,f) + T(r,g) + O(1) \leq 2T(r,g) + S(r,g).$$\\
Let $f_1 = (1/a) f,\;f_2 = h,\; f_3 = (-1/a) hg$. Then,
\begin{equation}\label{11}
\sum^{3}_{j=1} f_j \equiv 1.
\end{equation}
 Combining (\ref{b}), (\ref{c}) and (\ref{d}), we get
$$\sum^{3}_{j=1}\left(\overline{N}(r,f_j) + \overline{N}(r, \frac{1}{f_j})\right) = S(r,g).$$
 Clearly, $f_1, f_2 \;\text{and}\;f_3$ are linearly dependent and so there exist three constants $c_1, c_2$ and $c_3$ (atleast one of them is not zero) such that
\begin{equation} \label{12}
\sum^{3}_{j = 1} {c_j f_j} = 0
\end{equation}
If $c_1 = 0$, then from (\ref{12}) we see that $c_2\neq 0,\;c_3 \neq 0$, and
\begin{equation}\label{13}
f_3 = - \frac{c_2}{c_3}f_2.
\end{equation}
Substituting (\ref{13}) into (\ref{11}) gives
\begin{equation}\label{14}
f_1 + (1- \frac{c_2}{c_3})f_2 = 1.
\end{equation}
From (\ref{13}) and (\ref{14}), we get
$$ T(r, f_3) = T(r, f_1) + O(1)$$
 and thus
\begin{equation}\label{15}
T(r) = T(r, f_1) + O(1)
\end{equation}
where $T(r) ={\max} _{1 \leq j \leq 3}\{T(r, f_j)\}$.

Since $f_1$ is not a constant, it follows from (\ref{14}) that $1 - c_2/c_3 \neq 0$. From (\ref{14}), (\ref{15}) and Lemma \ref{0}, we deduce that
$$T(r) < \overline N\left(r, \frac{1}{f_1}\right) + \overline N\left(r, \frac{1}{f_2}\right) + \overline N(r, f_1) + S(r) = S(r),$$
where $S(r) = o(T(r))$, which is a contradiction and so $c_1 \neq 0$, and then (\ref{12}) gives
\begin{equation}\label{16}
f_1 = - \frac{c_2}{c_1}f_2 - \frac{c_3}{c_1} f_3.
\end{equation}
Now, from (\ref{11}) and (\ref{16}), we get
\begin{equation}\label{17}
\left(1 - \frac{c_2}{c_1}\right)f_2 + \left(1 - \frac{c_3}{c_1}\right)f_3 = 1.
\end{equation}
We consider the following three cases:\\

$Case\;1$: $1 - c_2/c_1 \neq 0$ and $1 - c_3/c_1 \neq 0.$

In this case,(\ref{16}) and (\ref{17}) gives
\begin{equation}\label{18}
f_1 = \frac{c_2 - c_3}{c_1 - c_2}f_3 - \frac{c_2}{c_1 - c_2}.
\end{equation}
From (\ref{17}) and (\ref{18}), we have
$$T(r, f_2) = T(r, f_1) + O(1)$$
 and hence
\begin{equation}\label{19}
T(r) = T(r, f_1) + O(1).
\end{equation}
Applying Lemma \ref{0} to (\ref{17}) and using (\ref{19}), we obtain
$$T(r) < \overline N\left(r, \frac{1}{f_2}\right) + \overline N\left(r, \frac{1}{f_3}\right) + \overline N(r, f_2) + S(r) = S(r),$$
which is a contradiction.\\

$Case\;2$: $1 - c_2/c_1 = 0$.

From (\ref{17}), we have $1 - c_3/c_1 \neq 0$, and

\begin{equation}\label{20}
f_3 = \frac{c_1}{c_1 - c_3}.
\end{equation}
Since $1 - c_2/c_1 = 0$, we obtain $c_1 = c_2$. Thus from (\ref{16}) and (\ref{20}), we obtain
\begin{equation} \label{21}
f_1 + f_2 = -\frac{c_3}{c_1 - c_3}.
\end{equation}
If $c_3 \neq 0$, then by applying Lemma \ref{0} to (\ref{21}), we obtain
$$T(r) < \overline N\left(r, \frac{1}{f_1}\right) + \overline N\left(r, \frac{1}{f_2}\right) + \overline N(r, f_1) + S(r) = S(r),$$
which is a contradiction. Hence $c_3 = 0$ and so from (\ref{20}), it follows that $f_3 \equiv 1$.\\

$Case\;3$: $1 - c_3/c_1 = 0$.

From (\ref{17}), we have $1 - c_2/c_1 \neq 0$, and

\begin{equation}\label{31}
f_2 = \frac{c_1}{c_1 - c_2}.
\end{equation}
Since $1 - c_3/c_1 = 0$, we obtain $c_1 = c_3$. Thus from (\ref{16}) and (\ref{31}), we obtain
\begin{equation} \label{32}
f_1 + f_3 = -\frac{c_2}{c_1 - c_1}.
\end{equation}
If $c_2 \neq 0$, then by applying Lemma \ref{0} to (\ref{32}), we obtain
$$T(r) < \overline N\left(r, \frac{1}{f_1}\right) + \overline N\left(r, \frac{1}{f_3}\right) + \overline N(r, f_1) + S(r) = S(r),$$
which is a contradiction. Hence $c_2 = 0$ and so from (\ref{31}), it follows that $f_2 \equiv 1$.\\

Thus if $f_2 \equiv 1$, then by (\ref{b}), we get, $f\equiv g$. If $f_3 \equiv 1$, then (\ref{b}) gives $f.g \equiv a^2.$ This proves $(i)$.

From Lemma \ref{22} , we see that $f^{(k)}, g^{(k)} \in A$. Using the conclusion of $(i)$, we get, either
  $$f^{(k)} \equiv g^{(k)}$$
  or
  $$f^{(k)}.g^{(k)} \equiv  a^2.$$
   If $f^{(k)} \equiv g^{(k)}$, then from Lemma \ref{33} ,we have $f \equiv g$. This completes the proof of $(ii)$.
\end{proof}

~~~~~~~~~~~~~~~~~~~~~~~~~~~~~~~~~~~~~~~~~~~~~~~~~~~~~~~~~~~~~~~~~~~~~~~~~~~~~~~~~~~~~~~~~~~~~~~~~~~~~~~~~~~~~~~~~$\Box$

\begin{proof}[Proof of Theorem 2.4]~~ Let the functions $F$ and $G$ be given by
$$F=\frac{f^{n+1}}{n+1}\; \text{and}\;G= \frac {g^{n+1}}{n+1}.$$
By hypothesis, $\overline E_{1)}(a, f^nf') = \overline E_{1)}(a, g^ng')$, therefore
 $$\overline E_{1)}(a, F')  = \overline E_{1)}(a, G').$$
Now
\begin{align*}
\overline{N}(r,F) + \overline{N}\left(r, \frac{1}{F}\right)
&  = \overline{N}\left(r,\frac{f^{n+1}}{n+1}\right) + \overline{N}\left(r, \frac{n+1}{f^{n+1}}\right)\\
& = \overline{N}(r,f) + \overline{N}(r, \frac{1}{f})\\
& = S(r,f)\\
& = S(r,F).
\end{align*}
Similarly by replacing $F$ by $G$ in above equation, we have
$$\overline{N}(r,G) + \overline{N}(r, \frac{1}{G}) = S(r,G).$$
Thus $F, G \in \mathcal{A}$ and so by the Theorem $2.1$, it follows that either
$$F'G' \equiv a^2 \;\text{or}\; F \equiv G.$$
Consider the case $F'G' \equiv a^2$, that is,
\begin{equation} \label{+}
 f^nf'g^ng' \equiv a^2.
 \end{equation}
Suppose that $z_1$ is a pole of $f$ of order $p$. Then $z_1$ is a zero of $g$ of order say $q$ and so from (\ref{+}), we find that
$$nq + q - 1 = np +p + 1. $$
That is, $(q-p)(n+1) = 2$, which is not possible as $n\geq 2$ and $p,q$ are positive integers. Thus $f$ and $g$ are entire functions and so from Lemma \ref{111}, we get $f(z)=c_1e^{cz}$ and $g(z)=c_2e^{-cz}$, where $c,c_1,c_2$ are constants satisfying $(c_1c_2)^{(n+1)}c^2 = -a^2$.

Next consider the case when $F \equiv G$. This gives
$$\frac{f^{n+1}}{n+1} = \frac{g^{n+1}}{n+1},$$
or
$$ f^{n+1} = g^{n+1}.$$
Hence $f = dg$ for some $(n+1)th$ root of unity $d$.
\end{proof}

~~~~~~~~~~~~~~~~~~~~~~~~~~~~~~~~~~~~~~~~~~~~~~~~~~~~~~~~~~~~~~~~~~~~~~~~~~~~~~~~~~~~~~~~~~~~~~~~~~~~~~~~~~~~~~~~~$\Box$

\begin{proof}[Proof of Theorem 2.6]~~
Let the functions $F$ and $G$ be given by
$$F = \frac{f^n(f^m)^{(k)}}{a}\;\text{ and } \;G = \frac{g^n(g^m)^{(k)}}{a}.$$
 Since $f^n(f^m)^{(k)}$ and $g^n(g^m)^{(k)}$ share $a$ CM, $F$ and $G$ share 1 CM.
Since (by Lemma \ref{55} and $T(r, a) = S(r, f)$),

\begin{align*}
N_2\left(r, \frac{1}{F}\right) + N_2(r, F)
& \leq N_2\left(r,\frac{1}{f^n(f^m)^{(k)}}\right) + N_2\left(r, f^n(f^m)^{(k)}\right) + S(r, f)\\
&\leq  N_2\left(r, \frac{1}{f^n}\right) + N_2\left(r, \frac{1}{(f^m)^{(k)}}\right) + 2 \overline N\left(r, f^n(f^m)^{(k)}\right) + S(r, f)\\
& \leq  2 \overline N\left(r, \frac{1}{f}\right) + N\left(r, \frac{1}{(f^m)^{(k)}}\right) + 2\overline N(r,f) + S(r, f)\\
& \leq  2 \overline N\left(r, \frac{1}{f}\right) + N\left(r, \frac{1}{f^m}\right) + k \overline N(r, f^m) + 2 \overline N(r, f) + S(r,f)\\
& = 2 \overline N\left(r, \frac{1}{f}\right) + m N\left(r, \frac{1}{f}\right) + k \overline N(r, f) + 2 \overline N(r, f) + S(r,f)\\
& = 2 \overline N\left(r, \frac{1}{f}\right) + mN\left(r, \frac{1}{f}\right) + (k + 2) \overline N(r, f) + S(r,f)\\
& \leq  2 T(r, f) + m T(r, f) + (k + 2)T(r,f) + S(r,f)\\
& =  (k + m + 4)T(r,f) + S(r, f),
\end{align*}
therefore,
\begin{equation} \label{5}
N_2\left(r, \frac{1}{F}\right) + N_2(r, F) \leq (k + m + 4)T(r,f) + S(r, f).
\end{equation}
On the similar lines we can write (\ref{5}) for the function $G$ as
\begin{equation} \label{6}
N_2\left(r, \frac{1}{G}\right) + N_2(r, G)  \leq (k + m + 4)T(r,g) + S(r, g).
\end{equation}
Since
\begin{align*}
nT(r, f) = T(r, f^n) & =  T\left(r,\frac{f^n(f^m)^{(k)}}{a}.\frac{a}{(f^m)^{(k)}} \right)\\
&\leq T(r, F) + T \left(r, \frac{1}{(f^m)^{(k)}} \right) + T(r, a) + S(r,f)\\
& \leq  T(r, F) + T \left(r, \frac{1}{(f^m)^{(k)}} \right) + S(r,f)\\
&\leq  T(r, F) + (k + 1)T \left(r, \frac{1}{f^m} \right) + S(r,f)\\
& =  T(r,F) + (km + m) T\left(r, \frac{1}{f}\right) + S(r,f),
\end{align*}
therefore
\begin{equation} \label{7}
(n - km -m )T(r, f) \leq T(r, F) + S(r,f).
\end{equation}
Similarly,
\begin{equation} \label{8}
(n - km -m )T(r, g) \leq T(r, G) + S(r,g).
\end{equation}
 Adding (\ref{7}) and  (\ref{8}), we get
 \begin{equation} \label{9}
 (n - km -m )\{T(r, f) + T(r, g)\} \leq \{T(r, F) + T(r, G)\} + S(r,f) + S(r, g).
 \end{equation}
Suppose that
 \begin{equation} \label{10}
 T(r, F) + T(r, G) \leq 2\{N_2\left(r, \frac{1}{F}\right) + N_2\left(r, \frac{1}{G}\right) + N_2(r,  F) + N_2(r, G\} + S(r,F) + S(r, G).
 \end{equation}
  holds. Then from (\ref{5}), (\ref{6}), (\ref{9}) and (\ref{10}), we have
  \begin{align*}
  (n-km-m)\{T(r,f) + T(r,g)\}  &\leq 2\{N_2\left(r, \frac{1}{F}\right) + N_2\left(r, \frac{1}{G}\right) + N_2(r,  F) + N_2(r, G)\} \\
   & \quad{} + S(r,f) + S(r,g).\\
  & \leq  2(k + m + 4)\{T(r,f) + T(r,g)\} + S(r,f) + S(r,g).\\
  & =  (2k + 2m + 8)\{T(r,f) + T(r,g)\} + S(r,f) + S(r,g),
  \end{align*}
  which implies that
  $$(n - km - 3m  -2k - 8)\{T(r,f) + T(r,g)\} \leq S(r,f) + S(r,g),$$
  a contradiction since $n > km + 3m + 2k + 8$, where $m > k-1$.

  Thus, by Lemma \ref{66}, it follows that either
  $$F.G \equiv 1 $$
  or
  $$F \equiv G.$$
   That is, either
   $$f^n(f^m)^{(k)}g^n(g^m)^{(k)} \equiv a^2$$
   or
   $$f^n(f^m)^{(k)} = g^n(g^m)^{(k)}.$$
\end{proof}

~~~~~~~~~~~~~~~~~~~~~~~~~~~~~~~~~~~~~~~~~~~~~~~~~~~~~~~~~~~~~~~~~~~~~~~~~~~~~~~~~~~~~~~~~~~~~~~~~~~~~~~~~~~~~~~~~$\Box$

\begin{proof}[Proof of Theorem 2.10]~~
Since
\begin{align*}
m\left(r, \frac{1}{f}\right) &= m\left(r, \frac{f^{(k)}}{f}.\frac{1}{f^{(k)}}\right)\\
& \leq  m\left(r, \frac{1}{f^{(k)}}\right) + m\left(r, \frac{f^{(k)}}{f}\right)\\
&= m\left(r, \frac{1}{f^{(k)}}\right) + S(r,f),
\end{align*}
therefore,
$$T(r,f) - N\left(r, \frac{1}{f}\right) \leq T(r,f^{(k)}) - N\left(r, \frac{1}{f^{(k)}}\right) + S(r,f),$$
and so
\begin{equation}\label{2}
N\left(r, \frac{1}{f^{(k)}}\right) \leq T(r,f^{(k)}) - T(r,f) + N\left(r, \frac{1}{f}\right) + S(r,f).
\end{equation}

Applying second fundamental theorem of Nevanlinna \cite[Theorem 2,5, p.47]{HAY} to the function $f^{(k)}$, we get
$$T(r, f^{(k)}) \leq  \overline N(r,f^{(k)}) + \overline N\left(r, \frac{1}{f^{(k)}}\right) + \overline N\left(r,\frac{1}{f^{(k)} - a}\right) + S(r,f^{(k)}).$$
That is,
\begin{equation} \label{3}
T(r, f^{(k)}) \leq  \overline N(r,f) + \overline N\left(r, \frac{1}{f^{(k)}}\right) + \overline N\left(r,\frac{1}{f^{(k)} - a}\right) + S(r,f).
\end{equation}

Since $N\left(r, 1/f''\right) = S(r,f)$, it follows from Lemma \ref{68} with $k = 1$ that
\begin{align*}
\overline N(r, f) + N_1\left(r,\frac{1}{f}\right) &\leq   N\left(r, \frac{1}{f''}\right) + \epsilon T(r,f) + S(r,f)\\
&= \epsilon T(r,f) + S(r,f).
\end{align*}

Thus, from (\ref{2}), (\ref{3}) and the fact that $f$ has finitely many simple zeros, we get

\begin{align*}
T(r, f)  & \leq  \overline N\left(r,\frac{1}{f^{(k)} - a}\right) + \overline N(r, f) + N\left(r,\frac{1}{f}\right) + S(r,f).\\
& \leq  N\left(r,\frac{1}{f^{(k)} - a}\right) + \overline N(r, f) + N\left(r,\frac{1}{f}\right) + S(r,f).\\
& =  N\left(r,\frac{1}{f^{(k)} - a}\right) + \overline N(r, f) + N_1\left(r,\frac{1}{f}\right) + \overline N\left(r,\frac{1}{f}\right) + S(r,f).\\
& \leq  N\left(r,\frac{1}{f^{(k)} - a}\right) +  \epsilon T(r,f) + \frac{1}{2} N\left(r,\frac{1}{f}\right) + S(r,f).\\
& \leq  N\left(r,\frac{1}{f^{(k)} - a}\right) + \epsilon T(r,f) + \frac{1}{2} T(r,f) + S(r,f).\\
& =  N\left(r,\frac{1}{f^{(k)} - a}\right) + (\frac{1}{2} + \epsilon) T(r,f) + S(r,f),
\end{align*}
which implies that
\begin{equation} \label{+}
(\frac{1}{2} - \epsilon)T(r, f) \leq N\left(r,\frac{1}{f^{(k)} - a}\right) + S(r,f).
\end{equation}
Taking $\epsilon = 1/4$ in (\ref{+}), we get
$$T(r,f) \leq 4 N\left(r,\frac{1}{f^{(k)} - a}\right) + S(r,f).$$
Hence $f^{(k)} - a$ has infinitely many zeros for $k = 1,\;2,\;3,\;\dotsc$.
\end{proof}

~~~~~~~~~~~~~~~~~~~~~~~~~~~~~~~~~~~~~~~~~~~~~~~~~~~~~~~~~~~~~~~~~~~~~~~~~~~~~~~~~~~~~~~~~~~~~~~~~~~~~~~~~~~~~~~~~$\Box$

\bibliographystyle{amsplain}

\end{document}